\documentclass[12pt]{article}
\usepackage{amssymb,amsmath}
\usepackage{hyperref}
\usepackage{graphicx}
\usepackage{color}

\begin{document}

\title{\Large\bf Sampling by incomplete cosine expansion of the sinc function: application to the Voigt/complex error function}

\author{
\normalsize\bf S. M. Abrarov\footnote{\scriptsize{Dept. Earth and Space Science and Engineering, York University, Toronto, Canada, M3J 1P3.}}\, and B. M. Quine$^{*}$\footnote{\scriptsize{Dept. Physics and Astronomy, York University, Toronto, Canada, M3J 1P3.}}}

\date{July 17, 2014}
\maketitle

\begin{abstract}
A new sampling methodology based on incomplete cosine expansion series is presented as an alternative to the traditional sinc function approach. Numerical integration shows that this methodology is efficient and practical. Applying the incomplete cosine expansion we obtain a rational approximation of the complex error function that with the same number of the summation terms provides an accuracy exceeding the Weideman\text{'}s approximation accuracy by several orders of the magnitude. Application of the expansion results in an integration consisting of elementary function terms only. Consequently, this approach can be advantageous for accurate and rapid computation.
\vspace{0.25cm}
\vspace{0.25cm}
\end{abstract}

\section {Introduction}
Sinc function is defined as \cite{Gearhart1990}
$$
{\rm{sinc}}\left( t \right) \equiv \left\{ 
\begin{aligned}
\frac{{\sin \left( t \right)}}{t}, \quad &t \ne 0\\
1, \quad &t = 0.
\end{aligned} 
\right.
$$
Sinc function is widely used in practical applications due to its remarkable sampling property \cite{Stenger2011}. In particular, using the sinc function we can approximate a function $f\left( t \right)$ by taking appropriate set of the sampling points $\left\{ {{t_n},f\left( {{t_n}} \right)} \right\}$ along some interval $\left[a,b\right]$ as \cite{Rybicki1989, Lether1998}:
\begin{equation}\label{eq_1}
f\left( t \right) = \sum\limits_{n =  - N}^N {f\left( {{t_n}} \right){\rm{sinc}}} \left( {\frac{\pi }{h}\left( {t - t{  _n}} \right)} \right) + \varepsilon \left( t \right),
\end{equation}
where $2N + 1$ is the total number of the sapling points, $h$ is the small parameter and $\varepsilon \left( t \right)$ is the error term. This remarkable sampling property of the sinc function can be used, for example, in numerical integration or differentiation \cite{Stenger2011}.

Although according to equation \eqref{eq_1} the sinc function can approximate quite accurately a function by sampling, its application is not always practically convenient. For example, suppose that $f\left( {t,x,y,z...} \right)$ is a multivariable function and suppose that it can be represented as a multiplication of two functions $f\left( {t,x,y,z...} \right) = g\left( t \right)h\left( {t,x,y,z...} \right)$. Then the expansion of the function $g\left( t \right)$ in terms of sinc function series \eqref{eq_1} leads to sum of integrals of kind
\begin{equation}\label{eq_2}
\int\limits_a^b {f\left( {t,x,y,z...} \right)} \,dt \approx \sum\limits_{n =  - N}^N {\int\limits_a^b {g\left( {{t_n}} \right){\rm{sinc}}} \left( {\frac{\pi }{h}\left( {t - {t_n}} \right)} \right)h\left( {t,x,y,z...} \right)dt}
\end{equation}
that may not be necessarily integrable or expressed in terms of the rational functions for efficient computation since the sinc function is not an elementary function.

It this work we propose a methodology of sampling based on the incomplete cosine expansion series of the sinc function that effectively resolves this problem and, as an example, we demonstrate its effectiveness in derivation of a new approximation of the Voigt/complex error function where only $16$ summation terms are sufficient to obtain accuracy better than ${10^{ - 12}}$ over almost all the complex domain of practical importance. This example illustrates how the incomplete cosine expansion of the sinc function can be applied as a useful and versatile technique in numerical integration.

\section{Sampling methodology}
\subsection{Incomplete cosine expansion}
Vi\'eta discovered that the sinc function can be represented as infinite product of the cosine functions \cite{Gearhart1990, Kac1959}:
\begin{equation}\label{eq_3}
{\rm{sinc}}\left( t \right) = \prod\limits_{m = 1}^\infty  {\cos \left( {\frac{t}{{{2^m}}}} \right)}.
\end{equation}
There is a product-to-sum identity \cite{Quine2013}
\begin{equation}\label{eq_4}
\prod\limits_{m = 1}^M {\cos \left( {\frac{t}{{{2^m}}}} \right)}  \equiv \frac{1}{{{2^{M - 1}}}}\sum\limits_{m = 1}^{{2^{M - 1}}} {\cos \left( {\frac{{2m - 1}}{{{2^M}}}t} \right)},
\end{equation}
where $M$ is an integer. Comparing identities \eqref{eq_3} and \eqref{eq_4} we can conclude, therefore, that
$$
{\rm{sinc}}\left( t \right) = \mathop {\lim }\limits_{M \to \infty } \frac{1}{{{2^{M - 1}}}}\sum\limits_{m = 1}^{{2^{M - 1}}} {\cos \left( {\frac{{2m - 1}}{{{2^M}}}t} \right)}.
$$
Consequently, at any finite integer $M$ we can refer to the right side of equation \eqref{eq_4} as incomplete cosine expansion of the sinc function. 

In sampling techniques involving the sinc function, it is quite common to replace its argument as $t \to {\pi} t/h$  (see for example \cite{Rybicki1989, Lether1998}). Following this convention, we rewrite the incomplete cosine expansion of the sinc function as the series
\begin{equation}\label{eq_5}
\frac{1}{{{2^{M - 1}}}}\sum\limits_{m = 1}^{{2^{M - 1}}} {\cos \left( {\frac{{\pi \left( {2m - 1} \right)}}{{{2^M}h}}t} \right)}.
\end{equation}

One useful property of this expansion occurs at infinitesimal $h$:
$$
\mathop {\lim }\limits_{h \to 0} \mathop {\lim }\limits_{M \to \infty } \frac{1}{{{2^{M - 1}}}}\sum\limits_{m = 1}^{{2^{M - 1}}} {\cos \left( {\frac{{\pi \left( {2m - 1} \right)}}{{{2^M}h}}t} \right)}  \equiv \delta \left( t \right)
$$
where $\delta \left( t \right)$ is the Kronecker delta, equal to $1$ at $t = 0$ and equal to $0$ at $t \ne 0$.  This interesting behaviour can be proved from the observation that:
$$
\mathop {\lim }\limits_{h \to 0} {\rm{sinc}}\left( {\frac{\pi }{h}t} \right) \equiv \delta \left( t \right).
$$

Since the integer $M$ in this application is implied to be finite, the incomplete cosine expansion \eqref{eq_5} is a periodic function with period $T$ that is dependent upon the values $M$ and $h$. A larger value of the integer $M$ increases the number of the summation terms in the incomplete cosine expansion \eqref{eq_5} that leads to an increase of the period $T$.

As the incomplete cosine expansion approximates the original sinc function ${\rm{sinc}}\left( {\pi t/h} \right)$ within the range $t \in \left[ { - T/4,T/4} \right]$ near the origin, we can write the approximation as
\begin{equation}\label{eq_6}
\frac{1}{{{2^{M - 1}}}}\sum\limits_{m = 1}^{{2^{M - 1}}} {\cos \left( {\frac{{\pi \left( {2m - 1} \right)}}{{{2^M}h}}t} \right)}  \approx {\rm{sinc}}\left( {\frac{\pi }{h}t} \right), \quad\quad - T/4 \le t \le T/4.
\end{equation}

Figure 1 shows the results calculated at $M = 5$ and $h = 0.25$. The period of the incomplete cosine expansion is $T = 16$. For comparison the original sinc function is also shown by blue line.

\begin{figure}[ht]
\begin{center}
\includegraphics[width=22pc]{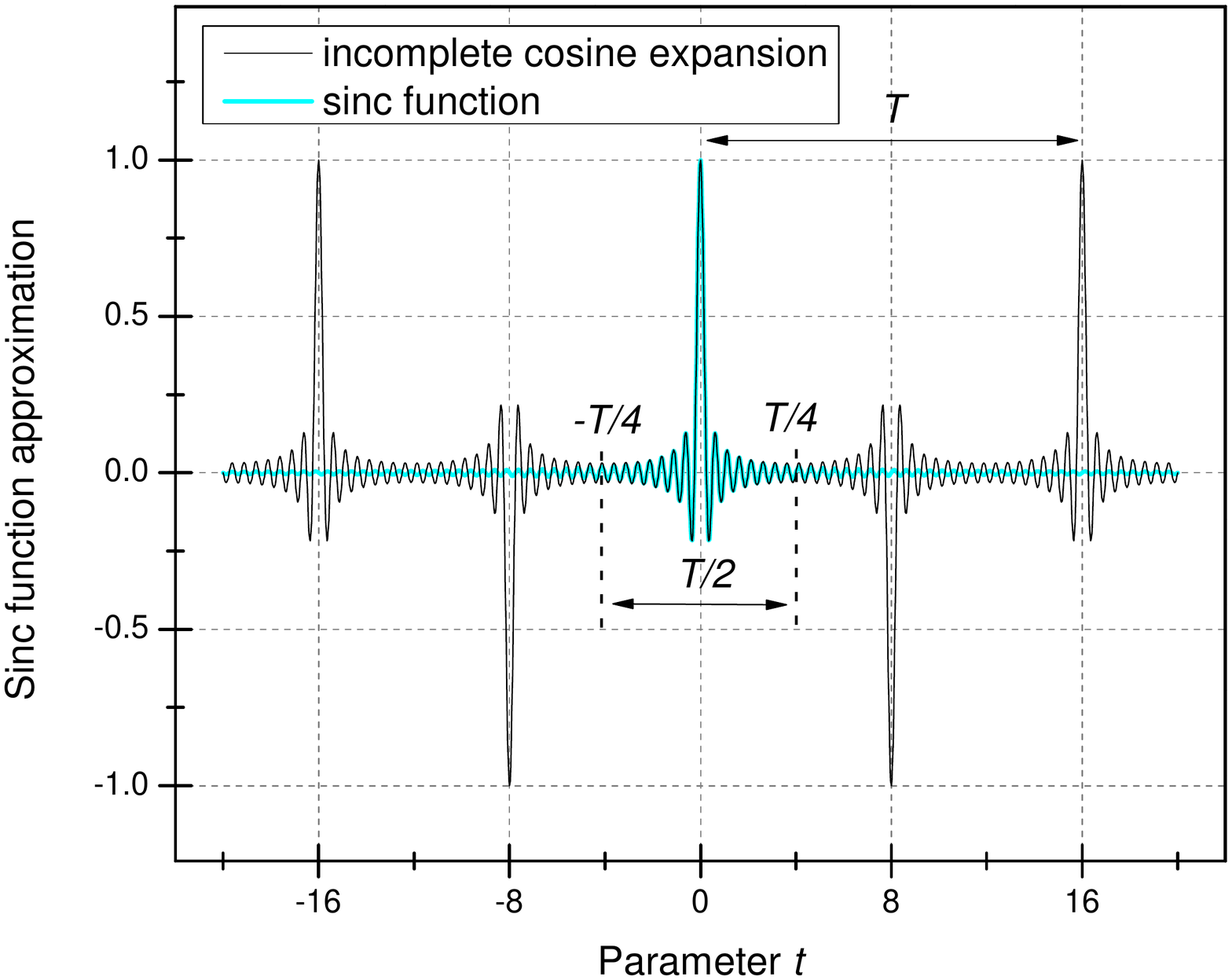}\hspace{2pc}%
\begin{minipage}[b]{24pc}
\vspace{0.3cm}
{\small{\sffamily {\bf{Fig. 1.}} The incomplete cosine expansion series at $M = 5$ and $h = 0.25$. The blue line shows the original sinc function ${\rm{sinc}}\left( {\pi t/h} \right)$.}}
\end{minipage}
\end{center}
\end{figure}

It should be noted that approximation \eqref{eq_6} is not the only incomplete cosine expansion of the sinc function. Examples of other incomplete cosine expansions are given in Appendix A.

\subsection{The Voigt/complex error function}
The efficiency of the incomplete cosine expansion can be demonstrated by application to the approximation of the complex error function defined as \cite{Schreier1992, Zaghloul2011}
$$
w\left( z \right) = {e^{ - {t^2}}}\left( {1 + \frac{{2i}}{{\sqrt \pi  }}\int\limits_0^z {{e^{{t^2}}}} dt} \right),
$$
where $z = x + iy$  is a complex argument.

Numerical approximations of the complex error function, also known as the Faddeeva function, find broad applications in many fields of Physics, Astronomy and Chemistry and its efficient algorithmic implementation stills remains a hot topic (see, for example, optimized algorithm in a recent work \cite{Karbach2014}). Mathematically, the real part of the complex error function $K\left( {x,y} \right) = {\mathop{\rm Re}\nolimits} \left[ {w\left( {x,y} \right)} \right]$, widely known as the Voigt function \cite{Schreier1992, Zaghloul2011}, represents a convolution of Cauchy (or Lorentz) and Gaussian (or Doppler) distributions. The Voigt function can be used to describe accurately spectral line broadening effects, for example, in atmospheric molecules \cite{Berk2013, Quine2002, Quine2013}, celestial bodies \cite{Emerson1996}, photo-luminescent materials \cite{Miyauchi2013} and so on. The simultaneous Lorentz and Doppler broadenings occur primarily due to the Heisenberg uncertainty principle, multiple collisions between particles and the velocity distribution of colliding particles.

The imaginary part of the complex error function $L\left( {x,y} \right) = {\mathop{\rm Im}\nolimits} \left[ {w\left( {x,y} \right)} \right]$ is also a useful function that can be used, for example, to describe the spectral behavior of refractive index in dispersive materials \cite{Balazs1969, Chan1986}.

Plasma Physics involving small-amplitude waves propagating through Maxwellian media is based on the plasma dispersion function as given by \cite{Fried1961, Swanson2003}
$$
Z\left( z \right) \equiv \frac{1}{{\sqrt \pi  }}\int\limits_{ - \infty }^\infty  {\frac{{{e^{ - {t^2}}}}}{{t - z}}} \,dt = i\sqrt \pi  w\left( z \right).
$$

The normal distribution function, applied in statistics and probability theory, can be expressed in terms of complex error function as \cite{Weisstein2003}
$$
\Phi \left( z \right) \equiv \frac{1}{{\sqrt {2\pi } }}\int\limits_0^z {{e^{ - {x^2/2}}}} dx = \frac{1}{2}\left[ {1 - {e^{ - {z^2}/2}}w\left( {\frac{{iz}}{{\sqrt 2 }}} \right)} \right].
$$

The error function of complex argument can also be expressed in terms of complex error function as follows
$$
{\rm{erf}}\left( z \right) = 1 - {e^{ - {z^2}}}w\left( {iz} \right).
$$
For example, the error function can be applied in a sub-pixel interpolation technique to star-image processing \cite{Quine2007}.

Consequently, the complex error function can be used as a versatile tool in Applied Mathematics and Computational Physics as the plasma dispersion function $Z\left( z \right)$, the normal distribution function $\Phi \left( z \right)$ and the error function ${\rm{erf}}\left( z \right)$ are simply the reformulated variations of the complex error function $w\left( z \right)$. Other functions that represent reformulations of the complex error function are the Fresnel integral \cite{Abramowitz1972, McKenna1984}, the Dawson\text{'}s integral \cite{Abramowitz1972, McKenna1984} and Gordeyev\text{'}s integral \cite{Swanson2003, Mace2003}.

As we discussed above, the application of the sinc function may not be always efficient. One of the representations of the complex error function is given by (see equation (3) in \cite{Srivastava1992}, see also \cite{Abrarov2011} for more details)
\begin{equation}\label{eq_7}
w\left( {x,y} \right) = \frac{1}{{\sqrt \pi  }}\int\limits_0^\infty  {\exp \left( { - {t^2}/4} \right)\exp \left( { - yt} \right)\exp \left( {ixt} \right)dt}, \quad\quad y > 0.
\end{equation}
Making change of variable $t/2 \to t$ in the integral \eqref{eq_7} leads to \cite{Abrarov2014}
\begin{equation}\label{eq_8}
w\left( {x,y} \right) = \frac{2}{{\sqrt \pi  }}\int\limits_0^\infty  {{e^{ - {t^2}}}{e^{2\left( {ixt - yt} \right)}}dt}, \quad\quad y > 0.
\end{equation}
As we can see, the integrand in the integral \eqref{eq_8} is a multivariable function. Obviously, this integrand can be rearranged as
$f\left( {t,x,y} \right) = g\left( t \right)h\left( {t,x,y} \right) = {e^{ - {t^2}}}{e^{2\left( {ixt - yt} \right)}}$,
where $g\left( t \right) = {e^{ - {t^2}}}$, $h\left( {t,x,y} \right) = {e^{2\left( {ixt - yt} \right)}}$. Consequently, according to equation \eqref{eq_2} we can approximate equation \eqref{eq_8} now as
\begin{equation}\label{eq_9}
w\left( {x,y} \right) \approx \frac{2}{{\sqrt \pi  }}\sum\limits_{n =  - N}^N {\int\limits_0^\infty  {{e^{ - t_n^2}}{\rm{sinc}}\left( {\frac{\pi }{h}\left( {t - {t_n}} \right)} \right){e^{2\left( {ixt - yt} \right)}}dt} }, \quad\quad		y > 0.
\end{equation}

The most simplest way is to take equidistantly separated sampling points ${t_n} = nh$ with step between two adjacent sampling points equal to $h$. With this assumption the approximation \eqref{eq_9} yields
\begin{equation}\label{eq_10}
w\left( {x,y} \right) \approx \frac{2}{{\sqrt \pi  }}\sum\limits_{n =  - N}^N {\int\limits_0^\infty  {{e^{ - {n^2}{h^2}}}{\rm{sinc}}\left( {\frac{\pi }{h}\left( {t - nh} \right)} \right){e^{2\left( {ixt - yt} \right)}}dt} }, \quad\quad y > 0.
\end{equation}

Although each integral term in the right side of equation \eqref{eq_10} is integrable through hyperbolic arctangent function, this method is absolutely inappropriate for rapid computation since the resultant approximation of the complex error function is not a rational approximation. 

This problem can be effectively resolved by using the incomplete cosine expansion. Specifically, we can write
\scriptsize
\begin{equation}\label{eq_11}
{e^{ - {t^2}}} \approx \frac{1}{{{2^{M - 1}}}}\sum\limits_{m = 1}^{{2^{M - 1}}} {\sum\limits_{n =  - N}^N {{e^{ - t_n^2}}\cos \left( {\frac{{\pi \left( {2m - 1} \right)\left( {t - {t_n}} \right)}}{{{2^M}h}}} \right)} }, \quad\quad - T/4 \le t \le T/4.
\end{equation}
\normalsize

Since the incomplete cosine expansion is a periodic function with period $T$, the obtained approximation \eqref{eq_11} is also periodic as it can be seen in the Fig. 2. Along with the pick at the origin coinciding with original exponential function ${e^{ - {t^2}}}$ shown by blue line, we can also observe additional negative and positive peaks at $t = \left\{ {\pm T/2,\pm 3T/2,\pm 5T/2\,,\pm 7T/2\,...} \right\}$ and $t = \left\{ {\pm T,\pm 2T,\pm 3T,\pm 4T,\pm 5T\,...} \right\}$, respectively. Therefore, the exponential function ${e^{ - {t^2}}}$ can be reasonably approximated only within the domain $ - T/4 \le t \le T/4$.

\begin{figure}[ht]
\begin{center}
\includegraphics[width=22pc]{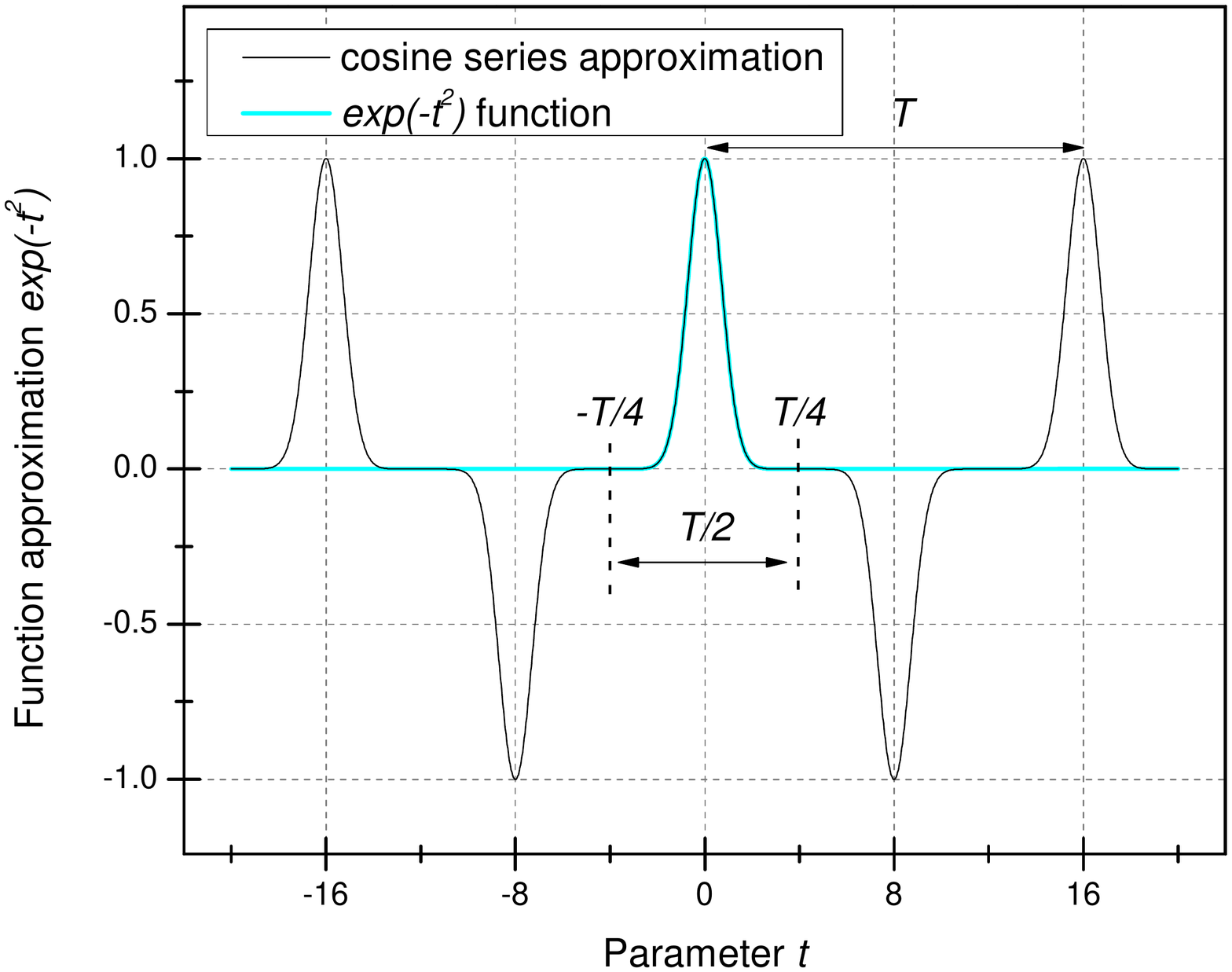}\hspace{2pc}%
\begin{minipage}[b]{24pc}
\vspace{0.3cm}
{\small{\sffamily {\bf{Fig. 2.}} Approximation of the function ${e^{ - {t^2}}}$by incomplete cosine expansion at $M = 5$ and $h = 0.25$. The blue line shows the original function ${e^{ - {t^2}}}$.}}
\end{minipage}
\end{center}
\end{figure}

Let us introduce a real constant $\varsigma  > 0$ and multiply both sides of equation \eqref{eq_11} by ${e^{ - \varsigma t}}$. If the constant $\varsigma $ is large enough, say equal or greater than $2$, the restriction $ - T/4 \le t \le T/4$ in equation \eqref{eq_11} can be effectively eliminated as ${e^{ - \varsigma t}}$ is a rapidly decreasing function that damps to zero all additional peaks as the argument $t$ increases. This tendency can be observed from Fig. 3 by comparing two curves at $\varsigma  = 0.1$ and $\varsigma  = 0.2$ shown by black and red colors, respectively. As the constant $\varsigma $ increases, all additional peaks are attenuated significantly. By $\varsigma  \mathbin{\lower.3ex\hbox{$\buildrel>\over
{\smash{\scriptstyle\sim}\vphantom{_x}}$}} 2$ they practically vanish and only the initial part of the approximated function ${e^{ - \varsigma t}}{e^{ - {t^2}}}$ remains unaffected. Thus, we can write
\begin{equation}\label{eq_12}
{e^{ - \varsigma t}}{e^{ - {t^2}}} \approx \frac{{{e^{ - \varsigma t}}}}{{{2^{M - 1}}}}\sum\limits_{m = 1}^{{2^{M - 1}}} {\sum\limits_{n =  - N}^N {{e^{ - t_n^2}}\cos \left( {\frac{{\pi \left( {2m - 1} \right)\left( {t - {t_n}} \right)}}{{{2^M}h}}} \right)} }, \quad\quad	\varsigma  \mathbin{\lower.3ex\hbox{$\buildrel>\over
{\smash{\scriptstyle\sim}\vphantom{_x}}$}} 2.
\end{equation}

\begin{figure}[ht]
\begin{center}
\includegraphics[width=22pc]{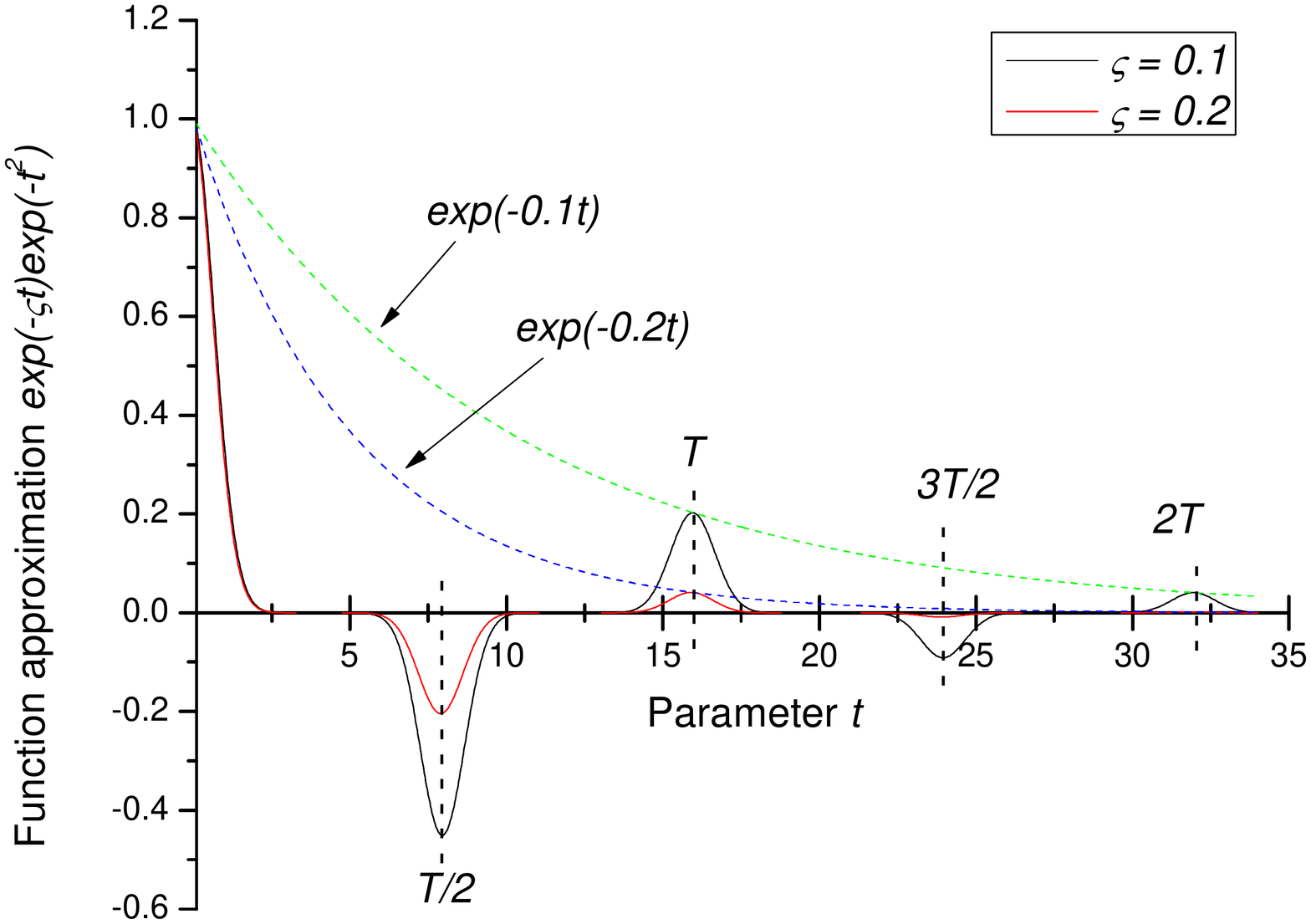}\hspace{2pc}%
\begin{minipage}[b]{24pc}
\vspace{0.3cm}
{\small{\sffamily {\bf{Fig. 3.}} Approximation of the function ${e^{ - {\varsigma ^2}}}{e^{ - {t^2}}}$by sampling with incomplete cosine expansion at $\varsigma  = 0.1$ and $\varsigma  = 0.2$ shown by black and red curves, respectively.}}
\end{minipage}
\end{center}
\end{figure}

From equation \eqref{eq_12} it follows that
\scriptsize
\begin{equation}\label{eq_13}
{e^{ - \varsigma t}}{e^{ - {{\left( {t - \varsigma /2} \right)}^2}}} \approx \frac{{{e^{ - \varsigma t}}}}{{{2^{M - 1}}}}\sum\limits_{m = 1}^{{2^{M - 1}}} {\sum\limits_{n =  - N}^N {{e^{ - t_n^2}}\cos \left( {\frac{{\pi \left( {2m - 1} \right)\left( {t - {t_n} - \varsigma /2} \right)}}{{{2^M}h}}} \right)} }, \quad\quad	\varsigma  \mathbin{\lower.3ex\hbox{$\buildrel>\over
{\smash{\scriptstyle\sim}\vphantom{_x}}$}} 2.
\end{equation}
\normalsize
As the peak of the exponential function ${e^{ - {{\left( {t - \varsigma /2} \right)}^2}}}$ is shifted to right from the origin by $\varsigma /2$, parameter $\varsigma $ can be regarded, therefore, as a shift constant.

Taking into account that
${e^{ - {t^2}}} = {e^{{\varsigma ^2}/4}}{e^{ - \varsigma t}}{e^{ - {{\left( {t - \varsigma /2} \right)}^2}}}$
and substituting approximation \eqref{eq_13} into \eqref{eq_8} leads to
\scriptsize
\begin{equation}\label{eq_14}
w\left( {x,y} \right) \approx \frac{{{e^{ {\varsigma ^2}/4}}}}{{{2^{M - 2}}\sqrt \pi  }}\sum\limits_{m = 1}^{{2^{M - 1}}} {\sum\limits_{n =  - N}^N {\int\limits_0^\infty  {{e^{ - t_n^2}}\cos \left( {\frac{{\pi \left( {2m - 1} \right)\left( {t - {t_n} - \varsigma /2} \right)}}{{{2^M}h}}} \right){e^{ - \varsigma t}}{e^{2\left( {ixt - yt} \right)}}dt} } }.
\end{equation}
\normalsize
Once again, assuming equidistant separation of the sampling points ${t_n} = nh$,  we can rewrite approximation \eqref{eq_14} in form
\scriptsize
\begin{equation}\label{eq_15}
w\left( {x,y} \right) \approx \frac{{{e^{ {\varsigma ^2}/4}}}}{{{2^{M - 2}}\sqrt \pi  }}\sum\limits_{m = 1}^{{2^{M - 1}}} {\sum\limits_{n =  - N}^N {\int\limits_0^\infty  {{e^{ - {n^2}{h^2}}}\cos \left( {\frac{{\pi \left( {2m - 1} \right)\left( {t - nh - \varsigma /2} \right)}}{{{2^M}h}}} \right){e^{ - \varsigma t}}{e^{2\left( {ixt - yt} \right)}}dt} } }.
\end{equation}
\normalsize
Each integral term on the right side of approximation \eqref{eq_15} is now integrable and can be expressed as a rational function. Consequently, from equation \eqref{eq_15} it follows that
\begin{equation}\label{eq_16}
w\left( z \right) \approx \sum\limits_{m = 1}^{{2^{M - 1}}} {\sum\limits_{n =  - N}^N {\frac{{{\alpha _{m,n}} + \left( {z + i\varsigma /2} \right){\beta _{m,n}}}}{{C_m^2 - {{\left( {z + i\varsigma /2} \right)}^2}}}} }, \quad\quad {\mathop{\rm Im}\nolimits} \left[ z \right] > 0,
\end{equation}
where
$$
{C_m} = \frac{{\pi \left( {2m - 1} \right)}}{{{2^{M + 1}}h}},
$$
$$
{\alpha _{m,n}} = \frac{{\sqrt \pi  \left( {2m - 1} \right){e^{{\varsigma ^2}/4 - {n^2}{h^2}}}}}{{{2^{2M}}h}}\sin \left( {\frac{{\pi \left( {2m - 1} \right)\left( {nh + \varsigma /2} \right)}}{{{2^M}h}}} \right)
$$
and
$$
{\beta _{m,n}} =  - i\frac{{{e^{{\varsigma ^2}/4 - {n^2}{h^2}}}}}{{{2^{M - 1}}\sqrt \pi  }}\cos \left( {\frac{{\pi \left( {2m - 1} \right)\left( {nh + \varsigma /2} \right)}}{{{2^M}h}}} \right).
$$

Although the approximation \eqref{eq_16} is accurate, it is not rapid due to the double summation. However, this problem can be readily resolved by defining the following constants
$$
{A_m} \equiv \sum\limits_{n =  - N}^N {{\alpha _{m,n}}} 
$$
and
$$
{B_m} \equiv \sum\limits_{n =  - N}^N {{\beta _{m,n}}}.
$$
Thus, after some trivial rearrangement we obtain
\begin{equation}\label{eq_17}
w\left( z \right) \approx \sum\limits_{m = 1}^{{2^{M - 1}}} {\frac{{{A_m} + \left( {z + i\varsigma /2} \right){B_m}}}{{C_m^2 - {{\left( {z + i\varsigma /2} \right)}^2}}}}, \quad\quad {\mathop{\rm Im}\nolimits} \left[ z \right] > 0,
\end{equation}
where
$$
{A_m} = \frac{{\sqrt \pi  \left( {2m - 1} \right)}}{{{2^{2M}}h}}\sum\limits_{n =  - N}^N {{e^{{\varsigma ^2}/4 - {n^2}{h^2}}}\sin \left( {\frac{{\pi \left( {2m - 1} \right)\left( {nh + \varsigma /2} \right)}}{{{2^M}h}}} \right)}
$$
and
$$
{B_m} =  - \frac{i}{{{2^{M - 1}}\sqrt \pi  }}\sum\limits_{n =  - N}^N {{e^{{\varsigma ^2}/4 - {n^2}{h^2}}}\cos \left( {\frac{{\pi \left( {2m - 1} \right)\left( {nh + \varsigma /2} \right)}}{{{2^M}h}}} \right)}.
$$

Approximation \eqref{eq_17} is a rational function since it contains no trigonometric or exponential functions dependent on the input parameters $x$ and $y$. Furthermore, all coefficients ${A_m}$, ${B_m}$ and ${C_m}$ involved in computation are independent of the input parameters $x$ and $y$. Therefore, the application of the approximation \eqref{eq_17} is advantageous for rapid computation.

A further advantage is that the approximation \eqref{eq_17} is readily integrable with respect to parameter $x$. As a result, it can be used for the spectrally integrated (or frequency integrated) Voigt function \cite{Bruggemann1992, Quine2013}.

It should be noted that for program coding the approximation \eqref{eq_17} can be conveniently rearranged in form
$$
w\left( z \right) \approx \psi \left( {z + i\varsigma /2} \right),
$$
where
$$
\psi \left( z \right) \equiv \sum\limits_{m = 1}^{{2^{M - 1}}} {\frac{{{A_m} + z{B_m}}}{{C_m^2 - {z^2}}}}
$$
is defined in order to deal with simplified $\psi \left( z \right)$-function.

To estimate the accuracy of this approximation \eqref{eq_17}, we define the relative errors for its real and imaginary parts as
$$
{\Delta _{{\mathop{\rm Re}\nolimits} }} = \left| {\frac{{{\mathop{\rm Re}\nolimits} \left[ {w\left( {x,y} \right)} \right] - {\mathop{\rm Re}\nolimits} \left[ {{w_{ref.}}\left( {x,y} \right)} \right]}}{{{\mathop{\rm Re}\nolimits} \left[ {{w_{ref.}}\left( {x,y} \right)} \right]}}} \right|
$$
and
$$
{\Delta _{{\mathop{\rm Im}\nolimits} }} = \left| {\frac{{{\mathop{\rm Im}\nolimits} \left[ {w\left( {x,y} \right)} \right] - {\mathop{\rm Im}\nolimits} \left[ {{w_{ref.}}\left( {x,y} \right)} \right]}}{{{\mathop{\rm Im}\nolimits} \left[ {{w_{ref.}}\left( {x,y} \right)} \right]}}} \right|,
$$
where ${w_{ref.}}\left( {x,y} \right)$ is the reference. The highly accurate reference values can be generated by using, for example, the program codes Algorithm 680 \cite{Poppe1990a, Poppe1990b} or recently published Algorithm 916 \cite{Zaghloul2011}.

At sufficiently large value $\left| z \right|$, say at  $\left| {x + iy} \right| \mathbin{\lower.3ex\hbox{$\buildrel>\over
{\smash{\scriptstyle\sim}\vphantom{_x}}$}} 15$, the computation of the complex error function is not problematic and many approximations are available in scientific literature for rapid and accurate computation (see, for example, the continuous fractional approximation in \cite{Gautschi1970}, the Gauss--Hermit quadrature \cite{Humlicek1979} or rational approximation in \cite{Hui1978}). Since in quantitative spectroscopy the range for the input parameters required for computation is $0 < x < 40,000$ and ${10^{ - 4}} < y < {10^2}$ \cite{Wells1999, Quine2002}, the most difficult domain for computation is $0 < x < 15$ and ${10^{ - 4}} < y < 15$. The proposed approximation \eqref{eq_17} of the complex error function very accurately covers this domain.

Figure 4a shows the logarithm of the relative error ${\log _{10}}{\Delta _{{\mathop{\rm Re}\nolimits} }}$ for the real part of the complex error function approximation \eqref{eq_17} calculated at $h = 0.25$, $N = 23$, $M = 5$ and $\varsigma  = 2.75$. In the most of the domain, the accuracy in real part is better than ${10^{ - 12}}$. Although in the real part the accuracy deteriorates at $y < {10^{ - 4}}$ , it still remains high and in the region ${10^{ - 4}} < y < {10^{ - 6}}$ it is better than ${10^{ - 8}}$.

Figure 4b illustrates the logarithm of the relative error ${\log _{10}}{\Delta _{{\mathop{\rm Im}\nolimits} }}$ for the imaginary part of the complex error function approximation \eqref{eq_17} also computed at $h = 0.25$, $N = 23$, $M = 5$ and $\varsigma  = 2.75$. In the most of the domain, the accuracy in imaginary part is better than ${10^{ - 12}}$. The accuracy deteriorates at $x < {10^{ - 4}} \cup y < {10^{ - 4}}$, but it is still high and better than ${10^{ - 8}}$. It should be noted that in radiative transfer spectroscopy the accuracy better than ${10^{ - 8}}$ is more than enough for practical applications.

Consider the Weideman\text{'}s approximation (see equation (38-I) and corresponding Matlab code in \cite{Weideman1994})
\begin{equation}\label{eq_18}
w\left( z \right) \approx \frac{{{\pi ^{ - 1/2}}}}{{L - iz}} + \frac{2}{{{{\left( {L - iz} \right)}^2}}}\sum\limits_{n = 0}^{N - 1} {{\gamma _{n + 1}}{{\left( {\frac{{L + iz}}{{L - iz}}} \right)}^n}}, \quad\quad {\mathop{\rm Im}\nolimits} \left[ z \right] > 0,
\end{equation}
where $L = {2^{ - 1/4}}{N^{1/2}}$ and
$$
{\gamma _n} = \frac{L}{\pi }\int\limits_{ - \infty }^\infty  {{e^{ - {t^2}}}{{\left( {\frac{{L - iz}}{{L + iz}}} \right)}^n}dt} 
$$

\begin{figure}[!htbp]
\begin{center}
\includegraphics[width=22pc]{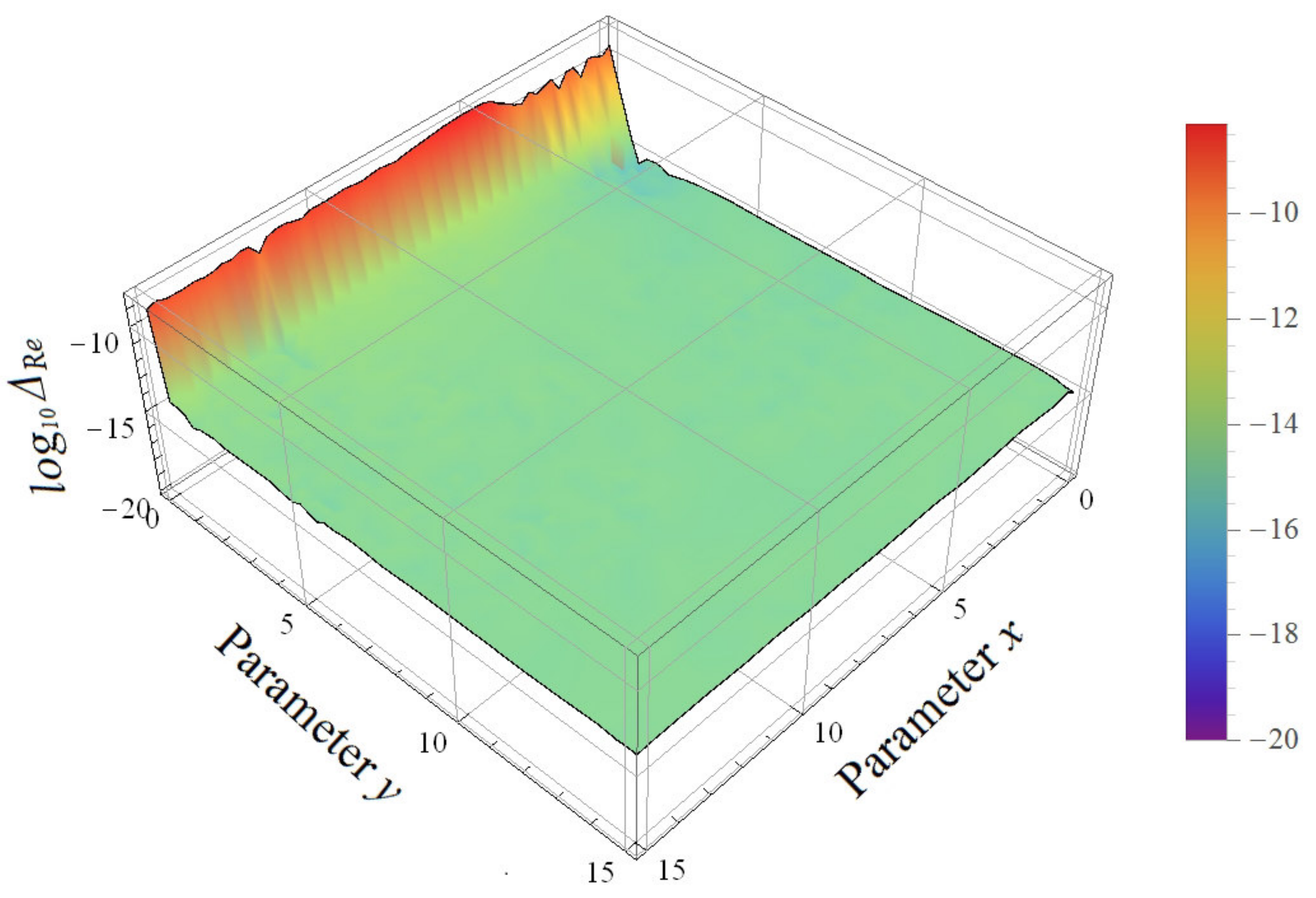}\hspace{2pc}%
\begin{minipage}[b]{24pc}
\vspace{0.25cm}
{\small{\sffamily {\bf{Fig. 4a.}} Logarithm of the relative error ${\log _{10}}{\Delta _{{\mathop{\rm Re}\nolimits} }}$ for the real part of approximation \eqref{eq_17} in the range $0 < x < 15$ and ${10^{ - 6}} < y < 15$.}}
\end{minipage}
\end{center}
\vspace{0.5cm}
\begin{center}
\includegraphics[width=22pc]{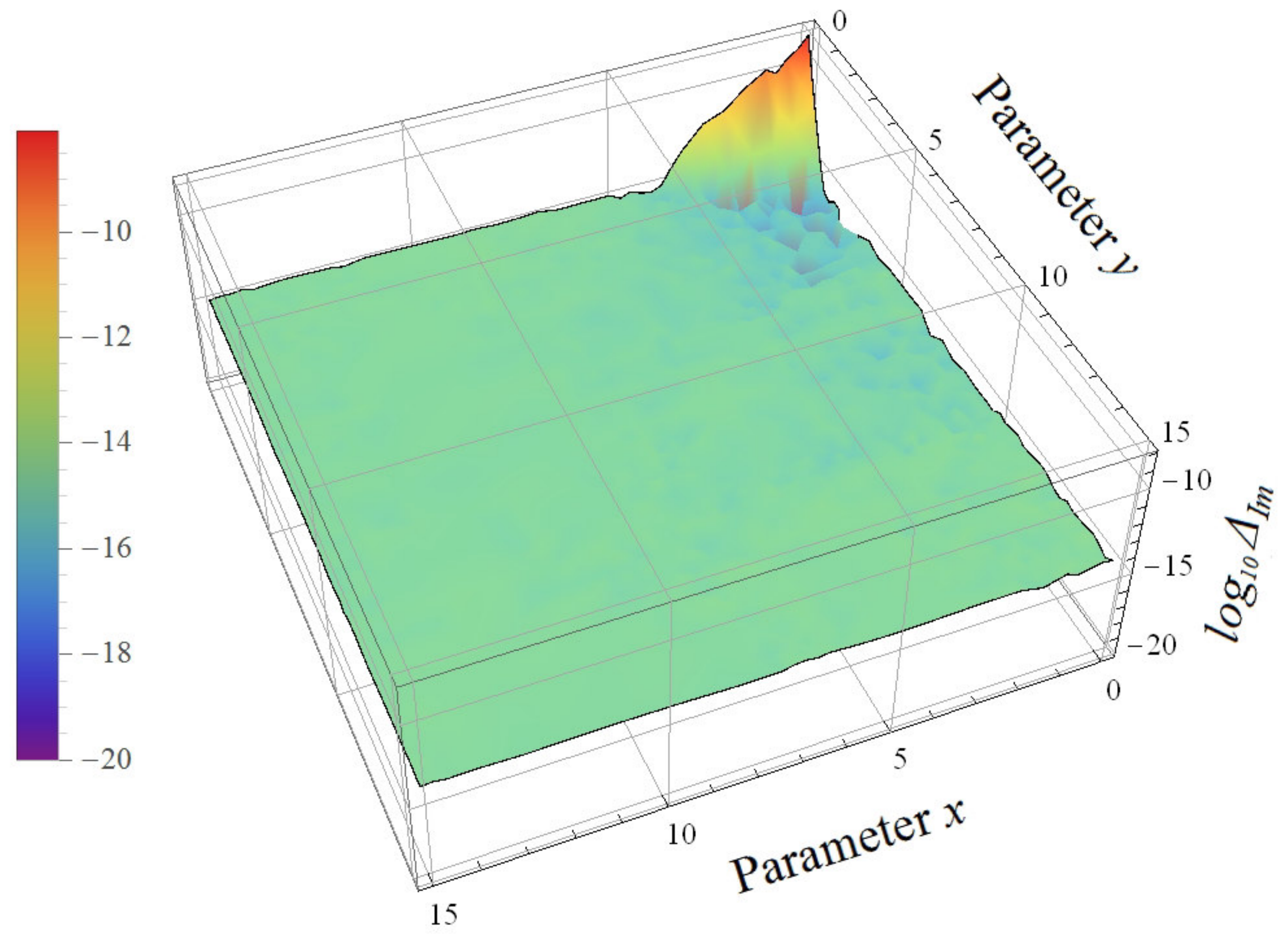}\hspace{2pc}%
\begin{minipage}[b]{24pc}
\vspace{0.25cm}
{\small{\sffamily {\bf{Fig. 4b.}} Logarithm of the relative error ${\log _{10}}{\Delta _{{\mathop{\rm Im}\nolimits} }}$for the imaginary part of approximation \eqref{eq_17} in the range $0 < x < 15$ and ${10^{ - 6}} < y < 15$.}}
\end{minipage}
\end{center}
\end{figure}

\noindent are the expansion coefficients that can be determined elegantly by FFT method. As Weideman\text{'}s approximation \eqref{eq_18} is a rational function, it is rapid in computation and, therefore, widely used in many applications. By default, Weideman applied the integer $N = 16$ in equation \eqref{eq_18}. 

Figure 5a shows the logarithm of the relative error ${\log _{10}}{\Delta _{{\mathop{\rm Re}\nolimits} }}$ for the real part of the complex error function approximation \eqref{eq_18}. In the most of the domain in real part the accuracy is about ${10^{ - 7}}$ and deteriorates at $y < {10^{ - 4}}$. Specifically, in the region ${10^{ - 4}} < y < {10^{ - 6}}$ the accuracy cannot be better than ${10^{ - 4}}$.

Figure 5b illustrates the logarithm of the relative error ${\log _{10}}{\Delta _{{\mathop{\rm Im}\nolimits} }}$ for the imaginary part of the complex error function approximation \eqref{eq_18}. The accuracy in the imaginary part is essentially better. Particularly, the accuracy in the most of the domain is about ${10^{ - 9}}$ and deteriorates to ${10^{ - 5}}$  in the region ${10^{ - 4}} < y < {10^{ - 6}}$.

\begin{figure}[!htbp]
\begin{center}
\includegraphics[width=22pc]{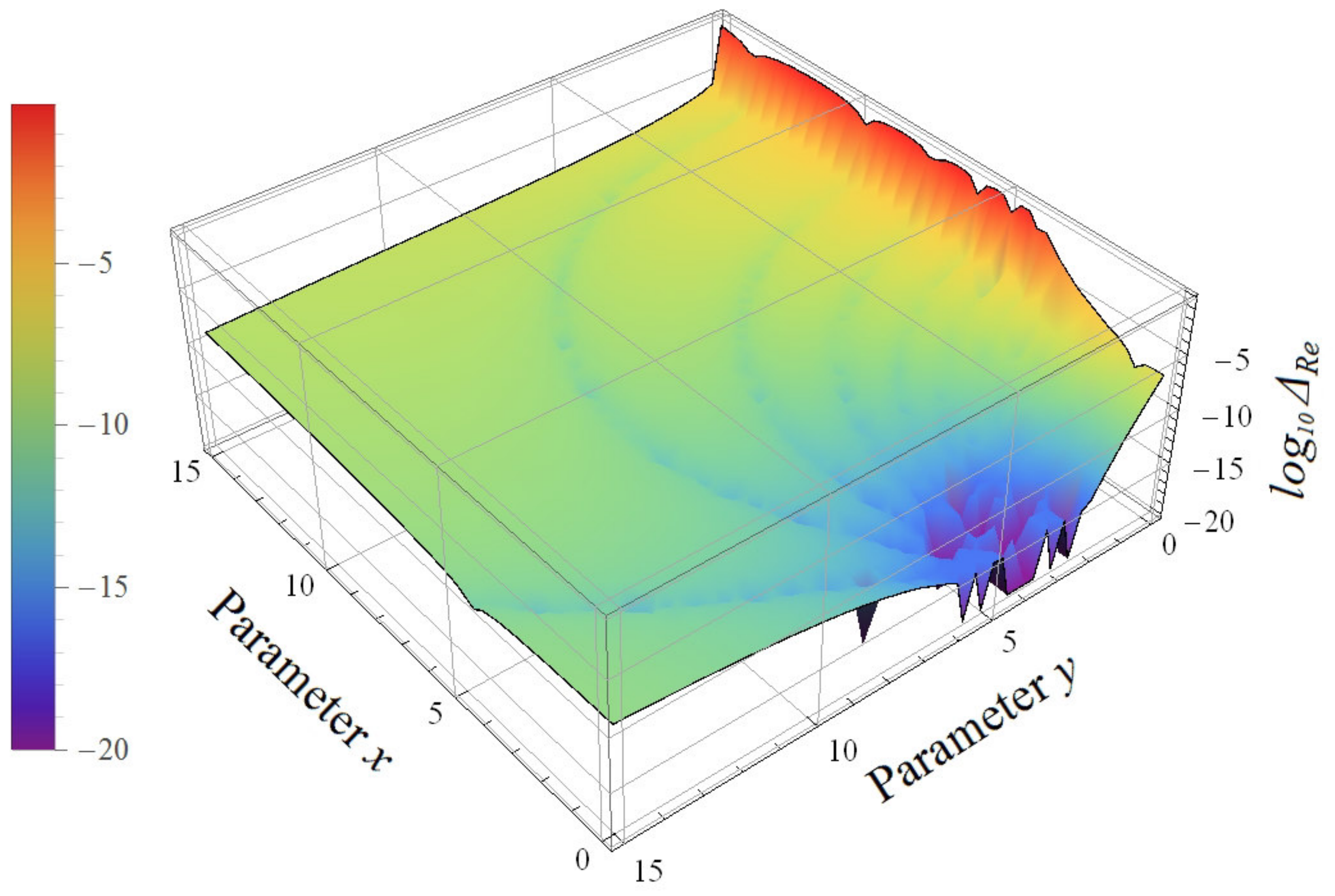}\hspace{2pc}%
\begin{minipage}[b]{24pc}
\vspace{0.25cm}
{\small{\sffamily {\bf{Fig. 5a.}} Logarithm of the relative error ${\log _{10}}{\Delta _{{\mathop{\rm Re}\nolimits} }}$ for the real part of Weideman\text{'}s approximation \eqref{eq_18} in the range $0 < x < 15$ and ${10^{ - 6}} < y < 15$.}}
\end{minipage}
\end{center}
\vspace{0.5cm}
\begin{center}
\includegraphics[width=22pc]{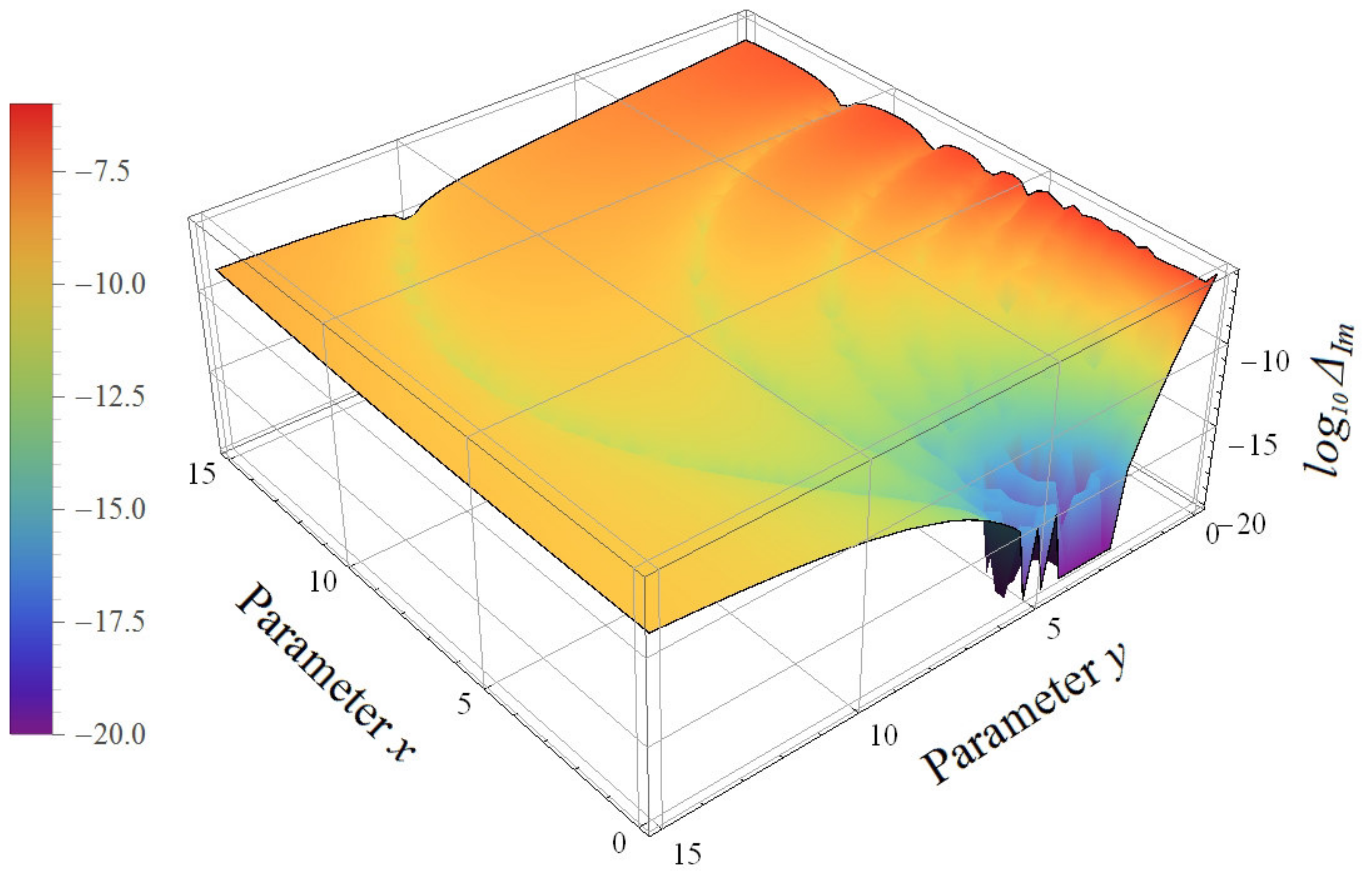}\hspace{2pc}%
\begin{minipage}[b]{24pc}
\vspace{0.25cm}
{\small{\sffamily {\bf{Fig. 5b.}} Logarithm of the relative error ${\log _{10}}{\Delta _{{\mathop{\rm Im}\nolimits} }}$ for the imaginary part of Weideman\text{'}s approximation \eqref{eq_18} in the range $0 < x < 15$ and ${10^{ - 6}} < y < 15$.}}
\end{minipage}
\end{center}
\end{figure}

As compared to the Weideman\text{'}s approximation \eqref{eq_18}, the approximation \eqref{eq_17} is more rapid in computation especially at extended size of the input array ${\bf{z}} = \left\{ {{z_1},{z_2},{z_3} \ldots } \right\}$. This is because the Weideman\text{'}s approximation \eqref{eq_18} is dependent on the power $n$. As a result of the exponentiations in polynomial terms, the Weideman\text{'}s approximation \eqref{eq_18} requires more computational resources when the size of the input array ${\bf{z}}$ exceeds several million elements.

Comparing Figs. 4a, 4b with 5a, 5b, one can see that while both approximations use 16 summation terms (${2^{5 - 1}} = 16$ in approximation \eqref{eq_17} and $N = 16$ in approximation \eqref{eq_18}), the complex error function approximation \eqref{eq_17} based on incomplete cosine expansion is more accurate than the Weideman\text{'}s approximation by several orders of the magnitude. These results illustrate that the application of the incomplete cosine expansion to a challenging problem can produce a highly accurate and efficient functional approximation.

\section{Conclusion}
We present a method of functional approximation based on incomplete cosine expansion series as a practical alternative to the traditional sinc function sampling. Our analysis shows that applying this approach to the complex error function we obtain an approximation that with same numbers of summation terms provides better accuracy than Weideman\text{'}s approximation by several orders of the magnitude. As application of the incomplete cosine expansion results in integration with only elementary functions, the approach is amenable to further mathematical flexibility and advantageous for accurate and rapid computation.

\section*{Acknowledgments}
This work is supported by Thoth Technology Inc., Discovery Channel Canada and York University.

\section* {Appendix A}
Since
$$
\cos \left( {\frac{{2m - 1}}{{{2^M}}}t} \right) = \cos \left( {\frac{{m - 1/2}}{{{2^{M - 1}}}}t} \right)
$$
the sinc function can also be represented as the limit
$$
{\text{sinc}}\left( t \right) = \mathop {\lim }\limits_{L \to \infty } \frac{1}{L}\sum\limits_{\ell  = 1}^L {\cos \left( {\frac{{\ell - 1/2}}{L}t} \right)}.
$$
This signifies that the upper value in summation may not be necessarily equal to $2^{M-1}$ and, in general, it can be taken as an arbitrary integer $L$:
\begin{equation}
\tag{A.1}\label{eq_A1}
\frac{1}{{{L}}}\sum\limits_{\ell = 1}^{{L}} {\cos \left( {\frac{{\pi \left( {\ell - 1/2} \right)}}{{{L}h}}t} \right)}  \approx {\rm{sinc}}\left( {\frac{\pi }{h}t} \right), \quad\quad - T_L/4 \le t \le T_L/4,
\end{equation}
where $T_L$ is the corresponding period that increases with increasing integer $L$. The equation \eqref{eq_6} is a specific case of a generalized equation \eqref{eq_A1} that occurs at $L = 2^{M-1}$.

Another expansion of the sinc function can be readily obtained by using the Poisson summation formula:
\begin{equation}
\tag{A.2}\label{eq_A2}
\sum\limits_{p =  - \infty }^{p = \infty } {{\text{sinc}}\left( {t + p{T_P}} \right)}  = \frac{\pi }{{{T_P}}}\sum\limits_{p =  -
\left\lfloor {{T_P}/\left( {2\pi } \right)} \right\rfloor }^{\left\lfloor {{T_P}/\left( {2\pi } \right)} \right\rfloor } {\cos
\left( {\frac{{2\pi p}}{{{T_P}}}t} \right)}
\end{equation}
where notation $\lfloor \; \rfloor$ denotes the floor function and  ${T_P}$ is the corresponding period. From this equation we obtain the following incomplete cosine expansion
$$
{\text{sinc}}\left( t \right) \approx \frac{\pi }{{{T_P}}}\sum\limits_{p =  - \left\lfloor {{T_P}/\left( {2\pi }
 \right)} \right\rfloor }^{\left\lfloor {{T_P}/\left( {2\pi } \right)}\right\rfloor } {\cos \left( {\frac{{2\pi p}}{{{T_P}}}t} \right)} , \quad\quad  - {T_P}/2 \leqslant t \leqslant {T_P}/2.
$$

It is interesting to note that from equation \eqref{eq_A2} it follows that
$$
\sum\limits_{p =  - \infty }^\infty  {{\text{sinc}}\left( {t + p{T_P}} \right)}  = \frac{\pi }{{{T_P}}}, \quad\quad 0 < \left| {{T_P}} \right|
 < 2\pi
$$
and at $T_P = 1$ we obtain an identity for the constant $\pi$:
$$
\sum\limits_{p =  - \infty }^\infty  {{\text{sinc}}\left( {t + p} \right)}  = \pi.
$$


\end{document}